\documentclass[12pt]{amsart}
\usepackage{tabmac}
\usepackage{amsthm,amssymb}
\usepackage{hyperref}
\usepackage{epsfig}
\usepackage{pstricks}
\usepackage{pst-plot}
\usepackage{pdflscape}
\usepackage{color}

\usepackage[margin=1.1in]{geometry}
\newtheorem{theorem}{Theorem}[section]

\newtheorem{lemma}[theorem]{Lemma}

\newtheorem{prop}[theorem]{Proposition}

\newtheorem{cor}[theorem]{Corollary}
\newtheorem{conjecture}[theorem]{Conjecture}

\theoremstyle{definition}
\newtheorem{remark}[theorem]{Remark}
\newtheorem{example}[theorem]{Example}

\def\Z{\mathbb{Z}}

\newcommand{\hako}{ 
\thicklines
\put(0,0){\line(0,1){0.8}}
\put(0,0){\line(1,0){0.8}}
\put(0.8,0){\line(0,1){0.8}}
\thinlines
}

\newcommand{\tama}{ 
\color[cmyk]{0,0,0,0.1}
\put(0.14,0.4){\circle*{0.6}}
\color{black}
\put(0.14,0.4){\circle{0.6}}
}

\newcommand{\kago}{ 
\put(0.05,0.3){\line(1,0){0.7}}
\qbezier(0.05,0.3)(0.05,0)(0.4,0)
\qbezier(0.4,0)(0.75,0)(0.75,0.3)
\qbezier(0.2,0.3)(0.2,0.2)(0.6,0.05)
\qbezier(0.4,0.3)(0.4,0.2)(0.7,0.15)
\qbezier(0.6,0.3)(0.6,0.2)(0.2,0.05)
\qbezier(0.4,0.3)(0.4,0.2)(0.1,0.15)
}

\author{Thomas Lam}
\email{tfylam@umich.edu}
\address{Department of Mathematics,
University of Michigan, 530 Church St., Ann Arbor, MI 48109 USA}
\author{Pavlo Pylyavskyy}
\email{ppylyavs@math.umn.edu}
\address{Department of Mathematics,
University of Minnesota, 206 Church St. SE, Minneapolis, MN 55455 USA}
\author{Reiho Sakamoto}
\email{reiho@rs.tus.ac.jp}
\address{Department of Physics,
Tokyo University of Science, Kagurazaka, Shinjukuku, Tokyo 162-8601 Japan}
\title{Box-basket-ball systems}
\begin{document}
\begin{abstract}
Using the whurl relation of the first two authors,
we define a new discrete solitonic system, which we call the box-basket-ball system,
generalizing the box-ball system of Takahashi and Satsuma.
In box-basket-ball systems balls may be put either into boxes or into baskets.
While boxes stay fixed, both balls and baskets get moved during time evolution.
Balls and baskets behave as fermionic and bosonic particles respectively.
We classify the solitons of this system, and study their scattering.
\end{abstract}
\maketitle

\section{Introduction}
In 1990 Takahashi--Satsuma \cite{TS,Tak} introduced a new
discrete soliton system called the box-ball system (BBS).
Their discovery is an outcome of the effort to find
a cellular automaton with solitonic behaviour,
which was a rather popular subject in the 1980s,
including the filter automata introduced by Park--Steiglitz--Thurston \cite{PST}.
A state of the BBS consists of an infinite sequence of boxes
where each box can accommodate at most one ball.
Then the time evolution rule of the BBS is described by a
simple combinatorial rule in terms of the box and ball interpretation of the system.

A striking feature of the BBS is that, despite its simple outlook,
it exhibits all characteristic properties of solitonic systems.
This remarkable property of the BBS is a sign of deep mathematical
structures behind the BBS.
In \cite{TTMS}, the authors realized the BBS as the ultra-discrete (or tropical)
limit of an ordinary soliton system, thereby proving the
integrability of the BBS.
Another line of the development is the generalization of the original BBS
by introducing extra degrees of freedom like spices of balls or
capacity of carriers \cite{Tak,TM}.
Such generalizations eventually culminated in the discovery
of the connection with Kashiwara's crystal bases theory \cite{Kas}
of quantum affine algebras found by many authors \cite{HHIKTT,FOY}.
Combining these two results, the BBS is now understood as both a
classical integrable system and a quantum integrable system.

Recently, yet another mathematical structure behind the BBS was revealed.
In \cite{KOSTY}, the inverse scattering formalism of the BBS is established
with the aid of the theory of rigged configurations.
Rigged configurations are certain combinatorial objects
originally introduced by Kerov--Kirillov--Reshetikhin \cite{KKR}
through their study of the Bethe ansatz analysis for quantum spin chains.
As an application, the initial value problem of the BBS
is solved \cite{KSY,S1} including all generalizations
of the BBS in \cite{HHIKTT,FOY}.

Nowadays, reversing the direction of these developments,
the BBS also helps to develop new mathematical theories.
For example, in \cite{S2}, the mysterious algorithm for the
bijection between rigged configurations and tensor products
of crystals of type $A^{(1)}_n$ is identified with the energy functions
of crystal bases theory via the time evolution operator of the BBS.
Another example is that the inverse scattering formalism of
the BBS gave a motivation for generalizing rigged configurations
to include not only highest weight elements
but also arbitrary elements in tensor products of crystals \cite{S3,DS}.
Such a generalization of the rigged configurations have also proven to be
useful on representation theory side.
Indeed, in \cite{OSS}, the rigged configuration gives an interesting
insight into the crystal structure of the Kirillov--Reshetikhin crystals
of type $D^{(1)}_n$ introduced by \cite{S4}.
These developments show that the BBS is not only a model in mathematical physics
but also gives a source for future mathematical theories.

\bigskip

The aim of the present article is to propose an entirely different
generalization of the BBS.
As discussed above, the BBS is related with the crystal bases theory.
To be more specific, the time evolution operators of the BBS is
described by the combinatorial $R$-matrix of the crystal bases theory.
In \cite{LP}, a generalization of the type $A^{(1)}_n$ combinatorial $R$-matrix,
which they called the whurl relation,
are introduced through the study of networks on a cylinder.

Let us explain the whurl relation in more detail.
The whurl relation is parametrized by vertical wires (pointing up)
and horizontal wires (pointing right or left) on a cylinder.
If all the horizontal wires point to the right,
the whurl relation reduces to the original combinatorial $R$-matrices
for the symmetric tensor product representation of type $A^{(1)}_n$.
On the other hand, if some of wires point to the left,
the whurl relation becomes a new map satisfying the Yang--Baxter relation.
Thus we can expect that the whurl relation will generate another class
of quantum integrable systems generalizing the BBS.

Among the whurl relation,
the simplest possible non-trivial extension is the three horizontal wires case,
where two of the wires point to the right and the other one points to the left
(see section \ref{sec:whurl}).
In the present article, we concentrate only on this three wire case
in order to clarify the most fundamental properties of the new system
without getting into the technical complexities concerning the general whurl relation.
In this way we stay close in spirit to the original Takahashi-Satsuma paper \cite{TS}.

The resulting new system, which we coin the {\bf box-basket-ball system} (BBBS),
has a remarkable novel property, namely,
it contains two entirely different kinds of particles whereas the BBS and
its generalizations have essentially one type of particles
(of various internal degrees of freedom).
In section \ref{sec:combinatorial}, we provide a combinatorial description
of the system in terms of boxes, balls and baskets,
generalizing the description for the BBS.
Here each box or basket can accommodate at most one ball,
whereas baskets can be put more than one on a box.
In this sense, we can regard balls as fermionic particles
and baskets as bosonic particles, together with mutual interaction
between balls and baskets.
In section \ref{sec:def_by_whurl}, we show that this combinatorial description
agrees with the definition in terms of the whurl relation.
As it turns out, the baskets originate from the wire pointing to the left
in the network on the cylinder.

Since the whurl relation satisfies the Yang--Baxter relation,
we see that our BBBS is a quantum integrable system
(Theorem \ref{quantum_integrability}).
However it is still a non-trivial problem to show that the BBBS
is indeed a soliton system.
The rest of the paper is devoted to show the solitonic property of the BBBS
as described in Theorem \ref{th:main}.
In order to achieve this goal, we first classify all possible
solitary waves that propagate without changing their shapes
if there is no scattering (Proposition \ref{prop:classification}).
A novel property is the presence of {\it slow solitons}.
This new kind of solitary waves are caused by a difference in the phase shifts
of the scatterings for two types of velocity one solitons as described in
Proposition \ref{prop:FF} and Corollary \ref{cor:SS}, respectively.
Finally we give a careful analysis of the general scatterings
in Proposition \ref{prop:FS} which forms the technical heart of the present paper.

We give a comment on related works.
In \cite{HI}, the authors constructed a supersymmetric extension of the BBS
by using the crystal bases for the quantum superalgebra introduced by \cite{BKK}.
Their system is different from ours since their extension is due to the addition
of fermionic particles to the original BBS whereas the BBBS is obtained by
adding bosonic particles to the BBS.

\bigskip

Let us say a few words about possible future directions of the present study.
Recall that in the original BBS case, the solitonic property of the system
is a sign of rich mathematical structures behind the model.
In our BBBS case, we also established the solitonic property of the system.
Thus it is not unreasonable to expect that the unknown symmetry behind the general
whurl relation might have deeper properties and
it will be worthy to clarify such underlying symmetry behind our system.

Another possibility we have in mind is for engineering purposes
related to transportation problems or traffic flow problems.
Since our system possesses obstacles described by baskets,
we hope that the integrability of our system might provide some tools
for detailed analysis of such problems. An alternative physical model would be 
waves in shallow water with some sand on the bottom. The sand represents the slow
bosonic particles, while the water represents the fast fermionic particles.  

\section{Discrete dynamical systems arising from whurl relations}

\subsection{Takahashi-Satsuma box-ball system}

In \cite{TS,Tak} Takahashi and Satsuma defined a discrete dynamical system exhibiting solitonic behaviour. 
We have a sequence of sites $\{S_i \mid i \in \Z\}$.
Each site $S_i$ contains one box which is either empty or has a ball inside. We denote the state of $S_i$ by a pair of numbers $(a_i, b_i)$, where $b_i$ is the number of balls in it
and $a_i$ is the number of extra balls that could fit in. Thus, each site is either in state $(1, 0)$, which we call the {\it {vacuum}} state, or in state $(0,1)$, which we call the {\it {ball}} state. The 
notation is set up so that one could also consider the case of more than one box at a site.  We shall however consider only the simple case of one box. 
We assume that $S_i = (1,0)$ for $i \ll 0$ or $i \gg 0$.

Time evolution of the system is as follows. The carrier travels from left ($i \ll 0$) to right ($i \gg 0$), having an infinite capacity to carry balls. She performs the following operations at each site $S_i$:
\begin{itemize}
 \item if $S_i=(0,1)$, she picks up this ball, changing the state into the vacuum $S_i=(1,0)$;
 \item if $S_i=(1,0)$, and she is carrying at least one ball, she drops one ball, changing the state of the site into $S_i=(0,1)$.
\end{itemize}
In other words, we consider each ball from left to right
(that is, starting from $S_{-\infty}$) and move the ball to the next available box. Each ball is moved exactly once. This completes the time evolution.
 
A {\it soliton} is a sequence of states which evolves at constant speed with no change in internal structure. 
A {\it basic soliton} is a soliton $A$ which cannot be decomposed into (non-trivial) solitons
$A'$ and $A''$ where $A'$ and $A''$ are separated by at least as many vacuum states as the speed of $A$.

\begin{theorem} \label{thm:bb} \cite{TS}
 \begin{enumerate}
  \item The basic solitons of box-ball system are strings of consecutive balls; the soliton of length $k$ has speed $k$.
  \item When a system consisting of a disjoint union of (suitably separated) basic solitons is allowed to scatter, the outcome is a sequence of basic solitons with the same set of lengths, arranged in non-decreasing order from left to right.
  \item Starting at any initial state, after a finite amount of time the system separates into basic solitons with non-decreasing length. 
 \end{enumerate}
\end{theorem}

One can also describe the action of the carrier as follows. Assume the carrier is in state $(a,b)$ and the site she is going through is in state $(c,d)$. Then the new state of the carrier and 
the site are given by:
\begin{align*}
a'&= a + \min(b, c) - \min(a, d)\\
b'&= b - \min(b, c) + \min(a, d)\\
c'&= c - \min(b, c) + \min(a, d)\\
d'&= d + \min(b, c) - \min(a, d)
\end{align*}

It is easy to check that when $a = \infty$ this gives the action of the carrier described above. Thus, one can assume that the carrier is originally in the state $(\infty, 0)$ and acts 
on each site according to the transformation above.

\subsection{Combinatorial $R$-matrix and whurl relations}
One can give a crystal base \cite{Kas} formulation of the box-ball systems \cite{HHIKTT,FOY} (see \cite{TM} for a combinatorial description).  The transformation 
describing the carrier action is a special case of the {\it {combinatorial $R$-matrix}} arising in the theory of crystals. In this context, states of the system are regarded as elements of
tensor products of the crystals
$\cdots \otimes b_{i-1} \otimes b_i \cdots \in
\cdots\otimes B_{k_{i-1}}\otimes B_{k_i}\otimes B_{k_{i+1}}\otimes\cdots$.

\begin{figure}[h!]
    \begin{center}
    \input{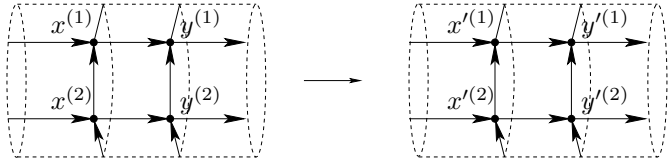}
    \end{center}
    \caption{An example of the whurl relation.}
    \label{fig:wire181}
\end{figure}

In \cite{LP} certain birational transformations constructed from networks on oriented surfaces were considered.  The tropicalization of one of these transformations turns out to be exactly the combinatorial $R$-matrix associated to the box-ball system.  More specifically, consider a cylinder with several parallel horizontal wires connecting
the components of the boundary, and several closed disjoint loops going around the cylinder.  All the loops are oriented in the same direction.  The orientations of the horizontal wires are allowed to vary:
each is oriented either from left to right or from right to left.  A parameter is associated to each vertex of the resulting network.
The case of two horizontal wires,
both oriented to the right, and two loops is shown in Figure \ref{fig:wire181}.

The loops around the cylinder are called {\it {whurls}}.
In \cite{LP} transformations of parameters of adjacent whurls were studied that preserve {\it {boundary 
measurements}} in the networks.  The {\it {whurl relation}} in the case shown in Figure \ref{fig:wire181} is given by
$$x'^{(1)} = y^{(1)} \frac{x^{(1)} + y^{(2)}}{y^{(1)} + x^{(2)}}; \;\;\; x'^{(2)} = y^{(2)} \frac{x^{(2)} + y^{(1)}}{y^{(2)} + x^{(1)}};$$
$$y'^{(1)} = x^{(1)} \frac{x^{(2)} + y^{(1)}}{y^{(2)} + x^{(1)}}; \;\;\; y'^{(2)} = x^{(2)} \frac{x^{(1)} + y^{(2)}}{y^{(1)} + x^{(2)}}.$$
Note that with $a = y^{(1)}, b=y^{(2)}, c =x^{(1)}, d= x^{(2)}$, this transformation recovers the piecewise-linear description of the box-ball system carrier, under the tropicalization (also called ultradiscretization) $(+,\times) \mapsto (\min,+)$.  We denote the whurl transformation by $R$.  

\begin{theorem} \cite[Proposition 11.8]{LP}
When all $n$ horizontal wires have the same orientation,
$R$ coincides with the birational version of the {\it {combinatorial $R$-matrix}}
of Kirillov-Reshetikhin crystals for symmetric powers of the standard representation
of $U_q'(A^{(1)}_n)$. 
\end{theorem}

In what follows we shall use the following property of $R$.  Suppose we are given three adjacent whurls labeled $1,2,3$ from left to right.  Let $R_{ij}$ be the whurl relation acting on the whurls labeled $i$ and $j$.
\begin{theorem}\label{YBE}
\cite[Theorem 6.6]{LP}
 The whurl relation (and thus, its tropicalization) satisfies the Yang-Baxter relation 
$$(R_{12} \otimes 1) \circ (1 \otimes R_{23}) \circ (R_{12} \otimes 1) =
(1 \otimes R_{23}) \circ (R_{12} \otimes 1) \circ (1 \otimes R_{23}).$$
\end{theorem}
%

\subsection{The mixed wires case and a box-basket-ball system} \label{sec:whurl}
The box-ball systems corresponding to whurl relations with all horizontal wires oriented in the same direction
have been extensively studied in literatures.
In particular, there are generalizations to systems with balls of many different colors, or to systems with boxes with higher capacity.  As suggested by the model of networks on surfaces in \cite{LP},
one may consider a more general setting of having wires in both directions.
The action of the carrier in such systems would be given 
by tropicalization of the corresponding whurl relations. We make the following general conjecture.

\begin{conjecture}
 The discrete dynamical systems arising from general whurl relations exhibit solitonic behaviour. 
\end{conjecture}

\begin{figure}[h!]
    \begin{center}
    \input{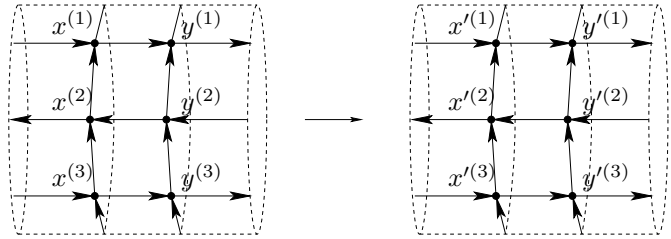}
    \end{center}
    \caption{A whurl relation with mixed directions of wires.}
    \label{fig:wire18}
\end{figure}

In this paper we shall consider the case of three horizontal wires, where two are oriented to the right and one to the left, as shown in Figure \ref{fig:wire18}. 
The whurl relation in this case is given by
$$x'^{(1)}=y^{(1)} \frac{x^{(1)}x^{(2)}+x^{(1)}x^{(3)}+x^{(2)}y^{(3)}}{y^{(2)}x^{(3)}+y^{(1)}x^{(3)}+y^{(1)}x^{(2)}}
\;\;\; x'^{(2)}=y^{(2)} \frac{x^{(1)}x^{(2)}+x^{(1)}x^{(3)}+x^{(2)}y^{(3)}}{x^{(1)}y^{(2)}+y^{(1)}y^{(3)}+y^{(2)}y^{(3)}}$$
$$ x'^{(3)}=y^{(3)} \frac{y^{(1)}x^{(3)}+y^{(1)}x^{(2)}+y^{(2)}x^{(3)}}{x^{(1)}y^{(2)}+y^{(1)}y^{(3)}+y^{(2)}y^{(3)}}$$
$$y'^{(1)}=x^{(1)} \frac{y^{(2)}x^{(3)}+y^{(1)}x^{(3)}+y^{(1)}x^{(2)}}{x^{(1)}x^{(2)}+x^{(1)}x^{(3)}+x^{(2)}y^{(3)}}
\;\;\; y'^{(2)}=x^{(2)} \frac{x^{(1)}y^{(2)}+y^{(1)}y^{(3)}+y^{(2)}y^{(3)}}{x^{(1)}x^{(2)}+x^{(1)}x^{(3)}+x^{(2)}y^{(3)}}$$
$$  y'^{(3)}=x^{(3)} \frac{x^{(1)}y^{(2)}+y^{(1)}y^{(3)}+y^{(2)}y^{(3)}}{y^{(1)}x^{(3)}+y^{(1)}x^{(2)}+y^{(2)}x^{(3)}}$$

In the next section we shall build a discrete dynamical system based on the tropicalization of this transformation, that we call the {\it {box-basket-ball system}}. 

\section{Box-basket-ball system}
\subsection{Combinatorial description}\label{sec:combinatorial}

We have a sequence of sites $\{S_i \mid i \in \Z\}$.
Each site $S_i$ contains one box, a number $b_i$ of baskets, and a number $c_i$ of balls.
Each ball must occupy a box or a basket, so that $c_i \leq b_i +1$.
We define $a_i = b_i-c_i+1$, the number of extra balls that can fit.
To denote the state of $S_i$, we will often use the vector $(a_i,b_i,c_i)$.
The vacuum state is $(1,0,0)$ -- one empty box and no baskets or balls.
We assume that $S_i = (1,0,0)$ for $i \ll 0$ or $i \gg 0$.

Time evolution of the system is as follows.  
Before beginning, we assume that a ball is always placed in a box if that is possible
(before being placed in a basket).
This does not change the state $(a_i,b_i,c_i)$.
First, we move every empty basket to the right one step (that is, from $S_i$ to $S_{i+1}$).
Full baskets are not moved.  Second, we consider each ball from left to right
(that is, starting from $S_{-\infty}$) and move the ball to the next available box or basket.
Each ball is moved exactly once.  This completes the time evolution.

\subsection{Piecewise linear description via carrier}\label{sec:def_by_whurl}
The time evolution can also be achieved via a carrier $C$,
which we think of as starting from the left,
and initially carrying infinitely many boxes, no baskets, and no balls.  
The interaction $(C,S_i) \mapsto (S'_i,C')$ of the carrier
$C = (a = \infty,b,c)$ with a site $S_i=(d,e,f)$ is given by
\begin{align*}
d'&= d+ b+f-\min(e+c,d+c,d+b)\\
e'&= e+ b-\min(d,e)\\
f'&= \min(e+c,d+c,d+b)-\min(d,e)
\end{align*}
where $S'_i=(d',e',f')$.  Balls and baskets are preserved, so that the resulting carrier $C'$ is given by 
\begin{align*}
a'&= \infty\\
b'&= \min(d,e)\\
c'&= c+f+\min(d,e) - \min(e+c,d+c,d+b).
\end{align*}
It is convenient to express the interaction of the carrier with
a site as the vertex diagram:
\begin{center}
\unitlength 12pt
\begin{picture}(22,7)
\put(3.5,0){$S_i'$}
\put(3.5,5.9){$S_i$}
\put(0.6,2.9){$C$}
\put(6.6,2.9){$C'$}
\multiput(0,0)(13,0){2}{
\put(4,1){\line(0,1){4.5}}
\put(1.7,3.2){\line(1,0){4.5}}
}
\put(9,3){or}
\put(11.5,3){$(a,b,c)$}
\put(19.5,3){$(a',b',c')$}
\put(15.4,5.9){$(d,e,f)$}
\put(15.2,0){$(d',e',f')$}
\end{picture}
\end{center}
\begin{prop}
The action of the carrier can be described as follows:
\begin{itemize}
 \item she picks up all the {\emph {empty}} baskets at the site,
and drops the baskets she was carrying before;
 \item then she drops as many of the balls she is carrying as she can,
and picks up {\emph {all}} the balls that were there before she dropped new ones.
\end{itemize}
\end{prop}

\begin{proof}
Initially there are $\min(d,e)$ empty baskets at the site,
which implies the correctness of the rule for baskets. 
After the baskets are already taken care of, there are $d - \min(d,e)$ empty boxes
and $b$ empty baskets at the site, and the correctness of the rule follows from the identity 
$$\min(e+c, d+c, d+b) - \min(d,e) = \min(c, d+b - \min(d,e)).$$
\end{proof}

\begin{prop}
The two descriptions of the time evolution are equivalent.
\end{prop}
\begin{proof}
The stated piecewise-linear transformation can be factorized as the composition of
\begin{align*}
d''&= d- \min(d,e)+b\\
e''&= e- \min(d,e)+b\\
f''&= f
\end{align*}
and
\begin{align*}
d'&= d''+f''-\min(c,d'')\\
e'&= e''\\
f'&= \min(c,d'')
\end{align*}
Indeed, as remarked above $$\min(c, d'') = \min(e+c, d+c, d+b) - \min(d,e),$$
and one needs only to plug in the expression for
$d'', e'', f''$ into the formulas for $d', e', f'$.

The two components of the piecewise-linear time evolution can be carried out independently of each other.
That is, one can first let the carrier go through all sites doing only the first
transformation, and then let her go through all sites doing only the second transformation.
Indeed, the first transformation does not depend on and does not influence the resulting
number of balls she carries, while the second transformation
does not depend on and does not influence the number of baskets she carries.
However, viewed independently from each other,
it is clear that the two transformations accomplish exactly
the combinatorial description of the time evolution.
\end{proof}

\begin{theorem}
 The action of the carrier in the box-basket-ball system is a tropicalization of the whurl relation given in Section \ref{sec:whurl}.
\end{theorem}

\begin{proof}
 Direct computation, plugging in $x^{(1)} = a = \infty$ into the tropicalization of the whurl relation.
\end{proof}

\subsection{Projection to box-ball system}
The box-ball system naturally embeds into the box-basket-ball system, by considering states with no baskets.  It turns out that one can also {\it project} the box-basket-ball system onto the box-ball system as follows.
Take a state of box-basket-ball system at time $t$.
Assume we have a site in a state $S = (a, b, c)$, where $a+c-b=1$.
Turn this site into a sequence of $b+1$ sites in the usual box-ball system, and fill the first $c$ of them with balls.  Do this for every site of the original state, creating a state of the box-ball system.
We call the operation the {\it {unbasketing}} of the original state; if the original state does not contain baskets, then unbasketing does not do anything.  Note however that there is no canonical way to reverse the unbasketing procedure.

\begin{theorem}
 The time evolution commutes with unbasketing. In other words, applying first the time evolution in the box-basket-ball system and then unbasketing gives the same result as first applying
unbasketing and then applying the time evolution in the box-ball system.
\end{theorem}

\begin{proof}
 The action of the carrier in the box-basket-ball system can be equivalently described the following way:
\begin{itemize}
 \item consider the empty baskets standing at the end of the unbasketing of a site $S_i$ as the beginning of the unbasketing of $S_{i+1}$;
 \item move the balls one by one left to right to the first unoccupied position;
 \item rearrange boxes and baskets without moving balls so that we get the correct unbasketing of the new state.
\end{itemize}
If we ignore the partitioning into sites in box-basket-ball system, as well as stop distinguishing between baskets and boxes, we are left only with the second step which is exactly 
the time evolution of the box-ball system.
\end{proof}

\begin{example}
 Start with the state $$\ldots, (1,2,2), (2,4,3), (1,2,2), \ldots,$$ where the sites not shown are vacuum. Its unbasketing is $\ldots 1, 1, 0, 1, 1, 1, 0, 0, 1, 1, \ldots$ where $1$ denotes a ball and $0$ denotes a vacuum state, and the sites not shown are vacuum states.  If we apply the time evolution to the original state we obtain $$\ldots, (2,1,0), (3,3,1), (2,3,2), (0,1,2), (0,0,1), (0,0,1) \ldots,$$ with unbasketing 
$\ldots, 1, 0, 0, 0, 1, 1, 0, 0, 1, 1, 1, 1, \ldots$. It is easily checked that this unbasketing is obtained from the original one by the time evolution of box-ball system. 
\end{example}

\section{Solitons}
\subsection{Classification}
Denote 
\begin{align*}
V &= (1,0,0) \\
F &=(0,0,1) \\
B_a &=(a+1,a,0), \;\;\; a \geq 1 \\
U_a &=(a,a,1), \;\;\; a \geq 1.
\end{align*}
\begin{example}
Let us depict a box by
{\unitlength 15pt
\begin{picture}(1,1)
\put(0,-0.15){\hako}
\end{picture}},
a ball by
{\unitlength 15pt
\begin{picture}(1,1)
\put(0,-0.15){\tama}
\end{picture}},
and a basket by
{\unitlength 15pt
\begin{picture}(1,1)
\put(0,0){\kago}
\end{picture}}.
Then the above elements are depicted as follows:
\begin{center}
\unitlength 15pt
\begin{picture}(15,2)
\put(0,0.15){$V=$}
\put(2,0){\hako}
\put(3,0){,}
\put(4,0.15){$F=$}
\put(6,0){\hako}
\put(6,0){\tama}
\put(7,0){,}
\put(8,0.15){$B_2=$}
\put(10,0){\hako}
\multiput(10,1)(0,0.5){2}{\kago}
\put(11,0){,}
\put(12,0.15){$U_2=$}
\put(14,0){\hako}
\put(14,0){\tama}
\multiput(14,1)(0,0.5){2}{\kago}
\put(15,0){.}
\end{picture}
\end{center}
\end{example}

A {\it soliton} is a sequence of states which evolves at constant speed with no change in internal structure. 
A {\it basic soliton} is a soliton $A$ which cannot be decomposed into (non-trivial) solitons
$A'$ and $A''$ where $A'$ and $A''$ are separated by at least as many vacuum states as the speed of $A$.

\begin{prop}\label{prop:classification}
The basic solitons are
\begin{enumerate}
\item
$F_k:=\overbrace{FF\cdots F}^k$ of speed equal to the length,
which we call a fast soliton,
\item 
any string of $F,B,U$ which does not contain the consecutive subsequence $FF$ or $FU$, of speed 1,
which we call a slow soliton.
\end{enumerate}
\end{prop}

\begin{proof}
Any soliton would remain a soliton after projecting to box-ball system via unbasketing.
Note that no basket can move with speed greater than one. Thus, if our soliton projects
to a soliton of speed greater than one, it has to contain no baskets.
A state containing no baskets projects to itself in the box-ball system. Thus by Theorem \ref{thm:bb} any 
such soliton would be one of the $F_k$, $k \geq 2$.

Assume now that the soliton has speed one.
Then it projects to a collection of disjoint solitons consisting of one ball.
Any positioning of boxes with at most one ball at each site 
has to be a combination of $F,B,U$ (we can assume there are no $V$-s since we are interested in basic solitons of speed one).
Since the unbasketting should not contain two consecutive balls,
the sequence should avoid $FF$ or $FU$ subsequences. 

It remains to argue that any such string of $F,B,U$ is indeed a soliton of speed one.
This follows from the carrier description of box-basket-ball system. We argue that carrier always carries at most one
ball to the next site, and always carries all the baskets from the previous site.
Indeed, the only way it could fail is if she is carrying a ball into a site and this site has already one
ball in it. However, since subsequences $FF$ and $FU$ are avoided,
she will also be carrying at least one basket, and thus would have where to put the ball.
\end{proof}

\subsection{Slower time evolution}

Consider the time evolution $T_{\ell}$ with carrier
$u_\ell :=(\ell, 0, 0)$ with $1 < \ell<\infty$.
The interaction $(C,S_i) \mapsto (S'_i,C')$ of the carrier $C = (a,b,c)$ with a site $S_i=(d,e,f)$ is given by
\begin{align*}
d'&= d+ \min(a+b, a+c, b+f) -\min(e+c, d+c, d+b)\\
e'&= e+ \min(a+b, a+c, b+f) -\min(a+e, d+f, e+f)\\
f'&= f+ \min(e+c, d+c, d+b) -\min(a+e, d+f, e+f)
\end{align*}
where $S'_i=(d',e',f')$.
Balls and baskets are preserved, so that the resulting carrier $C'$ is given by 
\begin{align*}
a'&= a - \min(a+b, a+c, b+f) + \min(e+c, d+c, d+b)\\
b'&= b - \min(a+b, a+c, b+f) + \min(a+e,d+f,e+f)\\
c'&= c - \min(e+c, d+c, d+b) + \min(a+e, d+f, e+f)
\end{align*}
Note that the relations $a'-b'+c'=a-b+c$ and
$d'-e'+f'=d-e+f$ hold.
Thus the capacities of sites and carriers are preserved under the time evolutions.
This is exactly the tropicalization of whurl relations in Section \ref{sec:whurl}.

\begin{prop}\label{prop:slow}
 The slower time evolution $T_{\ell}$ of solitons described above is as follows. 
\begin{enumerate}
 \item the soliton $F_k$ moves with speed $\min (\ell,k)$;
 \item the solitons of speed one evolve exactly the same way as before.
\end{enumerate}
\end{prop}

\begin{proof}
The first part is known \cite{FOY}, since in that case box-basket-ball system is indistinguishable from the box-ball system.

For the second part, observe that while acting on a speed one soliton the carrier never has to carry more than one ball at a time.
Thus it does not matter what her capacity is as
it is at least $1$. 
\end{proof}

Related to quantum integrability, our system
possesses the following property.
\begin{theorem}\label{quantum_integrability}
The time evolution operators $\{T_\ell\}_{\ell\geq 1}$ commute with each other.
\end{theorem}
For the proof, we prepare the following lemma.
\begin{lemma}
Let us consider a path $S=S_1\otimes S_2\otimes\cdots\otimes S_{m+n}$
where $S_j=V$ for $m\leq \forall j\leq m+n$.
If $n$ is sufficiently large,
we have $u_l\otimes S\simeq S'\otimes u_l$ under the isomorphism $R$.
\end{lemma}
\begin{proof}
Let us consider
$$(a,b,c)\otimes (1,0,0)\simeq
(d',e',f')\otimes (a',b',c')$$
under the whurl relation.
Recall that we have $a+c=b+l$ since we start from $u_l$.
Thus we have $a+c>b$, while $a+b\geq b$.
{}From this we have
\begin{align*}
b'&=b-\min (a+b,a+c,b)+\min (a,1,0)
=b-b=0.
\end{align*}
Similarly, if $c>0$, we have
$$c'=c-\min (c,1+c,1+b)+\min (a,1,0)=c-\min (c,b+1)<c.$$
Since $c$ is finite, if we consider the above isomorphism repeatedly,
we will finally obtain $c'=0$.
In this case, we have $(a',b',c')=(l,0,0)$,
which proves the lemma.
\end{proof}

\begin{proof}[Proof of Theorem \ref{quantum_integrability}]
Since $S_i=V$  for all $|i|\gg 1$, we can choose $m,n$
such that $S=S_{m}\otimes S_{m+1}\otimes\cdots\otimes S_n$
satisfies the above lemma as well as $S_i=V$ for all
$i<m$ and $n<i$.
We shall show $T_kT_l(S)=T_lT_k(S)$.
Apply the Yang--Baxter relation (Theorem \ref{YBE}) repeatedly to
$u_k\otimes u_l\otimes S$ as follows:
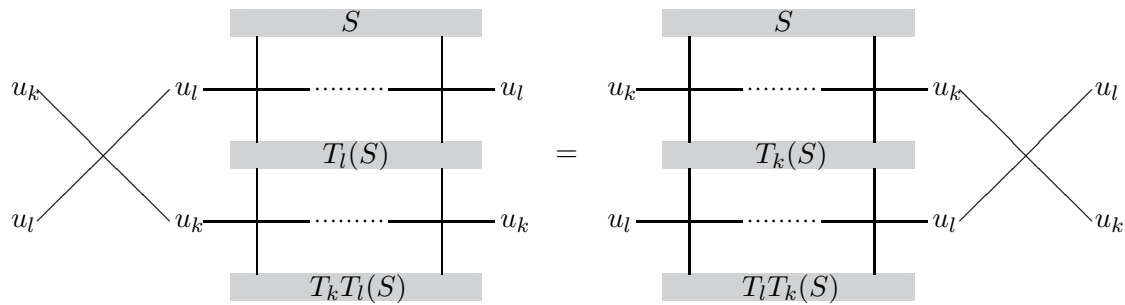
\begin{center}
{\small
\unitlength 10pt
\begin{picture}(34,11)
\color[cmyk]{0,0,0,0.2}
\put(4,0){\rule{95pt}{10pt}}
\put(4,5){\rule{95pt}{10pt}}
\put(4,10){\rule{95pt}{10pt}}
\color{black}
\put(1.7,3){\line(-1,1){5}}
\put(1.7,8){\line(-1,-1){5}}
\put(-4.3,2.8){$u_l$}
\put(-4.3,7.8){$u_k$}
\put(1.9,2.8){$u_k$}
\put(7,0.2){$T_kT_l(S)$}
\put(3,3){\line(1,0){4}}
\put(5,1){\line(0,1){4}}
\multiput(7.3,3)(0.3,0){9}{\circle*{0.1}}
\put(10,3){\line(1,0){4}}
\put(12,1){\line(0,1){4}}
\put(14.2,2.8){$u_k$}
\put(7.5,5.2){$T_l(S)$}
\put(1.9,7.8){$u_l$}
\put(3,8){\line(1,0){4}}
\put(5,6){\line(0,1){4}}
\multiput(7.3,8)(0.3,0){9}{\circle*{0.1}}
\put(10,8){\line(1,0){4}}
\put(12,6){\line(0,1){4}}
\put(14.2,7.8){$u_l$}
\put(8.2,10.2){$S$}
\put(16.3,5.3){$=$}
\put(16,0){
\color[cmyk]{0,0,0,0.2}
\put(4,0){\rule{95pt}{10pt}}
\put(4,5){\rule{95pt}{10pt}}
\put(4,10){\rule{95pt}{10pt}}
\color{black}
\put(1.9,2.8){$u_l$}
\put(7,0.2){$T_lT_k(S)$}
\put(3,3){\line(1,0){4}}
\put(5,1){\line(0,1){4}}
\multiput(7.3,3)(0.3,0){9}{\circle*{0.1}}
\put(10,3){\line(1,0){4}}
\put(12,1){\line(0,1){4}}
\put(14.2,2.8){$u_l$}
\put(7.5,5.2){$T_k(S)$}
\put(1.9,7.8){$u_k$}
\put(3,8){\line(1,0){4}}
\put(5,6){\line(0,1){4}}
\multiput(7.3,8)(0.3,0){9}{\circle*{0.1}}
\put(10,8){\line(1,0){4}}
\put(12,6){\line(0,1){4}}
\put(14.2,7.8){$u_k$}
\put(8.2,10.2){$S$}
}
\put(31.6,3){\line(1,1){5}}
\put(31.6,8){\line(1,-1){5}}
\put(36.7,2.8){$u_k$}
\put(36.7,7.8){$u_l$}
\end{picture}
}
\end{center}
Comparing both sides we obtain $T_kT_l(S)=T_lT_k(S)$.
\end{proof}
In the case of the original box-ball systems,
it is known that the family of mutually commutative time
evolutions gives complete information on the
action-angle variables of the dynamics \cite{S2}.

\section{Scattering}
\subsection{Some definitions}
We begin by introducing some notation that will be used in the rest of the paper.
The following two basic solitons we shall call {\it pure solitons}:
\begin{align*}
&F_k=\overbrace{FF\cdots F}^k,\\
&B_{a_1,a_2,\cdots, a_r}:=B_{a_1}B_{a_2}\cdots B_{a_r},
\end{align*}
and the remaining speed one basic solitons of Proposition
\ref{prop:classification} (2) are called {\it composite solitons}.

Let us introduce some notation for our speed 1 basic soliton.
The balls in the soliton we shall call {\it slow balls}.
If there is a ball in the tail (leftmost site) of a soliton,
then we call it the {\it initial slow ball}.
The condition that the soliton does not contain $FF$ or $FU$ is equivalent to
the condition that every non-initial slow ball has a basket
in the site immediately behind (towards the tail) of it.
Let us begin from the tail of the soliton and number all the baskets with $1,2,\ldots,b$.
Within the same site, we number from lower baskets to upper baskets.
If there is a non-initial slow ball in site $S_i$,
then we shall call the basket with the highest number in site $S_{i-1}$ {\it special}.
The $m$-th non-initial slow ball in the soliton
and the $m$-th special basket in the soliton are considered {\it paired}.

The balls in our speed $k$ soliton are called {\it fast balls}.

\subsection{Phase shift of scatterings}
Suppose $A$ and $A'$ are basic solitons,
such that after scattering one obtains again the solitons $A'$ and $A$
(but now in the reverse order): 
\begin{align*}
p=\cdots A\cdots\cdots A'\cdots\longrightarrow
T_\ell^N(p)=\cdots A'\cdots\cdots A\cdots
\end{align*}
for sufficiently large $N$.
Here, ``$\cdots A\cdots$" is the abbreviation for the state
$\cdots V\otimes V\otimes A\otimes V\otimes V\otimes\cdots$.
Let $v$ and $v'$ be the velocities of the free propagations of $A$ and $A'$
under the time evolution $T_\ell$.
Define the phase shifts of $A$ and $A'$ before the scattering to be 0.
Let the position of the leftmost letter of $A$ of $p$ be 0.
Then the position of the leftmost letter of $A$ of $T_\ell^N(p)$
can be expressed as $Nv+\delta$.
Similarly we define $\delta'$ for $A'$.
We call these $\delta$, $\delta'$ as phase shifts and describe
the above scattering as follows:
\begin{align*}
A[0]\ast A'[0]\longrightarrow
A'[\delta' ]\ast A[\delta ].
\end{align*}
If internal degree of freedom changes after the scattering,
we define the phase shift for each ball or basket in the similar way as above
by comparing with the corresponding elements before scattering.

\subsection{Scattering of fast solitons}
The following result for scattering in the box-ball system is well known.
\begin{prop}\label{prop:FF}
For any time evolution $T_\ell$ ($\ell\geq 2$),
the two body scatterings between the fast solitons is as follows ($m>n$):
\begin{align*}
F_m[0]\ast F_n[0]&\longrightarrow F_n[-2n]\ast F_m[2n]
\end{align*}
\end{prop}

\subsection{Scattering of fast soliton and slow soliton}
The general form of such scattering is as follows:
\begin{prop}\label{prop:FS}
For any time evolution $T_\ell$ ($\ell\geq 2$),
the two body scatterings between the fast soliton $F_m:\; m > 1$ and a basic slow soliton $A$ is given by:
\begin{align*}
F_m[0]\ast A&\longrightarrow A'\ast F_m[2b-a],
\end{align*}
where $b$ is the total number of balls in $A$, and $a$ is the total number of baskets in $A$, and $A'$ is obtained from $A$ as follows:
\begin{enumerate}
\item
All non-special baskets, and non-initial slow balls are slowed by 1.
\item
All special baskets are not phase-shifted.
\item
The initial slow ball (if any) is slowed by 2.
\end{enumerate}
\end{prop}

Note that because of (3), $A'$ may no longer be a basic soliton.

\begin{proof}
For $T_\infty$, this result is established in Section \ref{sec:FS}.
To show that the result of scattering using $T_\ell$ does not depend on $\ell$,
we use the commutativity of the time evolutions (Theorem \ref{quantum_integrability}).
If $p$ denotes the initial state then we have $(T_\ell)^M(p) = (T_\infty)^{-N}(T_\ell)^M (T_\infty)^N(p)$
where $N$ is chosen to be very large and $M$ is chosen even larger.
Suppose the scattering for $T_\infty$ is known,
so that $(T_\infty)^N(p)$ is the disjoint union of basic solitons as described in the proposition.
Then $(T_\ell)^M(T_\infty)^N(p)$ consists of the same set of solitions evolving freely at slow ($\leq \ell$) speed.
If $M \gg N$, then $(T_\infty)^{-N}(T_\ell)^M (T_\infty)^N(p)$ will still consist of
the same set of solitions with no (reverse) scattering.
Thus the scattering of $T_\ell$ and $T_\infty$ are the same.
\end{proof}

\begin{cor}\label{cor:SS}
Let $a$ be the total number of baskets in $B_{a_1,a_2,\cdots,a_r}$.
Then we have
$$F_m[0]\ast B_{a_1,a_2,\cdots, a_r}[0]\longrightarrow
B_{a_1,a_2,\cdots, a_r}[-1]\ast F_m[-a]$$
under any time evolution $T_\ell$ ($\ell\geq 2$).
\end{cor}

\begin{example}
Consider the scattering of $F_3$ with $B_1 U_3 F$.
Under $T_\infty$, it is displayed in Figure \ref{fig:T_infty}
\begin{figure}
\unitlength 15pt
\begin{picture}(20,3)
\put(0,0){\tama}
\put(1,0){\tama}
\put(2,0){\tama}
\put(6,0){\tama}
\put(7,0){\tama}
\multiput(0,0)(1,0){20}{\hako}
\put(5,1){\kago}
\multiput(6,1)(0,0.5){3}{\kago}
\end{picture}
\begin{picture}(20,3)
\put(3,0){\tama}
\put(4,0){\tama}
\put(5,0){\tama}
\put(7,0){\tama}
\put(8,0){\tama}
\multiput(0,0)(1,0){20}{\hako}
\put(6,1){\kago}
\multiput(7,1)(0,0.5){3}{\kago}
\end{picture}
\begin{picture}(20,3.5)
\put(6,0){\tama}
\put(7,0){\tama}
\put(8,1){\tama}
\put(8,0){\tama}
\put(9,0){\tama}
\multiput(0,0)(1,0){20}{\hako}
\put(7,1){\kago}
\put(8,1){\kago}
\multiput(8,2)(0,0.5){2}{\kago}
\end{picture}
\begin{picture}(20,3)
\put(8,0){\tama}
\put(9,0){\tama}
\put(9,1){\tama}
\put(10,0){\tama}
\put(11,0){\tama}
\multiput(0,0)(1,0){20}{\hako}
\multiput(8,1)(0,0.5){2}{\kago}
\multiput(9,1)(0,1){2}{\kago}
\end{picture}
\begin{picture}(20,3)
\put(9,0){\tama}
\put(10,0){\tama}
\put(12,0){\tama}
\put(13,0){\tama}
\put(14,0){\tama}
\multiput(0,0)(1,0){20}{\hako}
\multiput(9,1)(0,0.5){3}{\kago}
\put(10,1){\kago}
\end{picture}
\begin{picture}(20,3)
\put(10,0){\tama}
\put(11,0){\tama}
\put(15,0){\tama}
\put(16,0){\tama}
\put(17,0){\tama}
\multiput(0,0)(1,0){20}{\hako}
\multiput(10,1)(0,0.5){3}{\kago}
\put(11,1){\kago}
\end{picture}
\caption{Scattering of $F_3$ and $B_1U_3F$ under $T_\infty$.}
\label{fig:T_infty}
\end{figure}
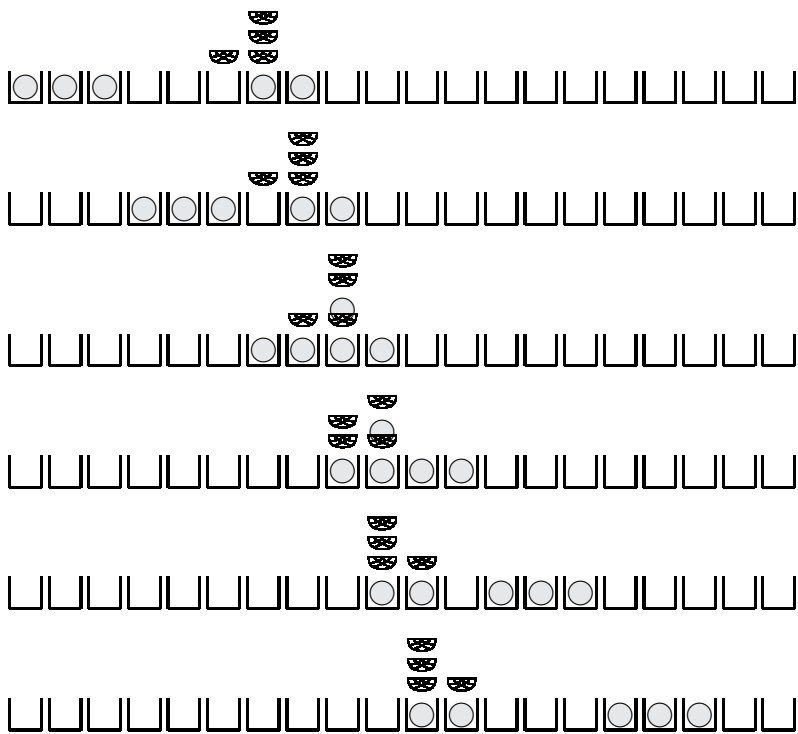
and under $T_2$ it is displayed in Figure \ref{fig:T_2}.
\begin{figure}
\unitlength 15pt
\begin{picture}(20,3)
\put(0,0){\tama}
\put(1,0){\tama}
\put(2,0){\tama}
\put(6,0){\tama}
\put(7,0){\tama}
\multiput(0,0)(1,0){20}{\hako}
\put(5,1){\kago}
\multiput(6,1)(0,0.5){3}{\kago}
\end{picture}
\begin{picture}(20,3)
\put(2,0){\tama}
\put(3,0){\tama}
\put(4,0){\tama}
\put(7,0){\tama}
\put(8,0){\tama}
\multiput(0,0)(1,0){20}{\hako}
\put(6,1){\kago}
\multiput(7,1)(0,0.5){3}{\kago}
\end{picture}
\begin{picture}(20,3)
\put(4,0){\tama}
\put(5,0){\tama}
\put(6,0){\tama}
\put(8,0){\tama}
\put(9,0){\tama}
\multiput(0,0)(1,0){20}{\hako}
\put(7,1){\kago}
\multiput(8,1)(0,0.5){3}{\kago}
\end{picture}
\begin{picture}(20,3)
\put(6,0){\tama}
\put(7,0){\tama}
\put(8,0){\tama}
\put(9,0){\tama}
\put(10,0){\tama}
\multiput(0,0)(1,0){20}{\hako}
\put(8,1){\kago}
\multiput(9,1)(0,0.5){3}{\kago}
\end{picture}
\begin{picture}(20,3.5)
\put(8,0){\tama}
\put(9,0){\tama}
\put(10,1){\tama}
\put(10,0){\tama}
\put(11,0){\tama}
\multiput(0,0)(1,0){20}{\hako}
\put(9,1){\kago}
\put(10,1){\kago}
\multiput(10,2)(0,0.5){2}{\kago}
\end{picture}
\begin{picture}(20,4)
\put(10,0){\tama}
\put(11,0){\tama}
\put(11,1){\tama}
\put(11,2){\tama}
\put(12,0){\tama}
\multiput(0,0)(1,0){20}{\hako}
\put(10,1){\kago}
\put(11,3){\kago}
\multiput(11,1)(0,1){2}{\kago}
\end{picture}
\begin{picture}(20,3)
\put(11,0){\tama}
\put(12,0){\tama}
\put(12,1){\tama}
\put(13,0){\tama}
\put(14,0){\tama}
\multiput(0,0)(1,0){20}{\hako}
\multiput(11,1)(0,0.5){2}{\kago}
\multiput(12,1)(0,1){2}{\kago}
\end{picture}
\begin{picture}(20,3)
\put(12,0){\tama}
\put(13,0){\tama}
\put(14,0){\tama}
\put(15,0){\tama}
\put(16,0){\tama}
\multiput(0,0)(1,0){20}{\hako}
\multiput(12,1)(0,0.5){3}{\kago}
\put(13,1){\kago}
\end{picture}
\begin{picture}(20,3)
\put(13,0){\tama}
\put(14,0){\tama}
\put(16,0){\tama}
\put(17,0){\tama}
\put(18,0){\tama}
\multiput(0,0)(1,0){20}{\hako}
\multiput(13,1)(0,0.5){3}{\kago}
\put(14,1){\kago}
\end{picture}
\caption{Scattering of $F_3$ and $B_1U_3F$ under $T_2$.}
\label{fig:T_2}
\end{figure}
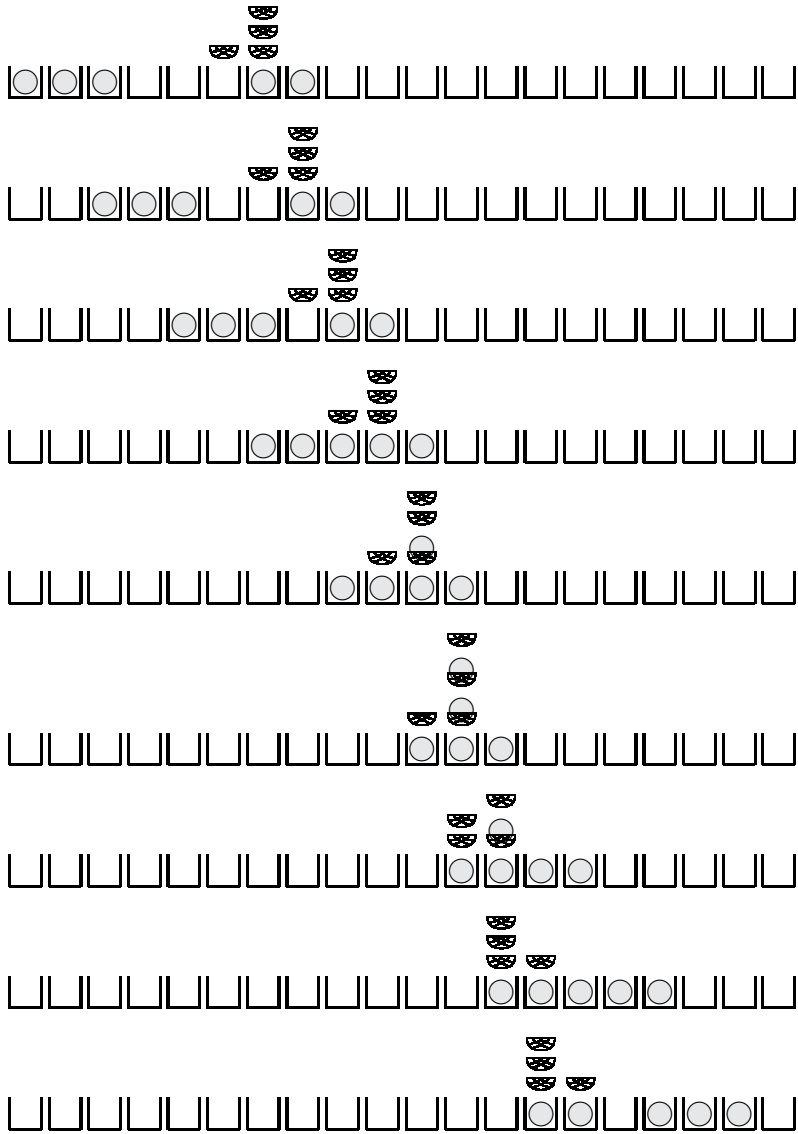
\end{example}

Any basic soliton of speed one can be uniquely cut into the chunks of the following kinds:
\begin{align*}
&B_a U_{b_1} \ldots U_{b_r} Z,\quad
B_a U_{b_1} \ldots U_{b_r} F Z,\quad
V U_{b_1} \ldots U_{b_r} Z,\quad
V U_{b_1} \ldots U_{b_r} F Z,\\
&V F Z,\quad
B_{a_1} \ldots B_{a_r} Z.
\end{align*}
Here $Z$ denotes either $B$ or $V$, though we regard that $Z$ is {\it {not}} a member of the chunk before $Z$. 


We can decompose composite basic solitons into disjoint union of pure solitons
by multiple scatterings with $F_k$'s.  (By a disjoint union  we mean that the basic solitions are separated by at least as many vacuum sites as the speed of the soliton.)

\begin{prop}\label{prop:composite}
For any time evolution $T_\ell$ ($\ell\geq 2$),
consider multiple scattering of composite solitons with
$F_k$ ($k\geq 2$) for sufficiently many times.
Then the chunks described above decompose into disjoint union
of pure solitons as follows:
\begin{align*}
B_a U_{b_1} \ldots U_{b_r}&
\longmapsto F^{\sqcup r}\sqcup B_{a+b_1+\cdots b_r}\\
B_a U_{b_1} \ldots U_{b_r} F&
\longmapsto F^{\sqcup r+1}\sqcup B_{a+b_1+\cdots b_r}\\
V U_{b_1} \ldots U_{b_r}&
\longmapsto F^{\sqcup r}\sqcup B_{b_1+\cdots b_r}\\
V U_{b_1} \ldots U_{b_r} F&
\longmapsto F^{\sqcup r+1}\sqcup B_{b_1+\cdots b_r}\\
V F&
\longmapsto F\\
B_{a_1,\cdots,a_r} &
\longmapsto B_{a_1,\cdots,a_r}
\end{align*}
Here $B_0 = V$ and $F^{\sqcup r}:=\overbrace{F\sqcup F\sqcup\cdots\sqcup F}^r$.
\end{prop}
\begin{proof}
The result essentially follows from Proposition \ref{prop:FS}, except that we must note that scattering of $F_k$ with $F \sqcup F \sqcup \cdots \sqcup F \sqcup A$ gives $F \sqcup F \sqcup \cdots \sqcup F \sqcup A'$.  In other words, the presence of the $F$'s do not affect the scattering of $F_k$ with $A$.  This is easy to establish by induction on the number of $F$'s. 
\end{proof}

Given a disjoint union of basic solitons, we define the number of {\it ball solitons} and number of {\it basket solitons} as follows.
For pure solitons, we consider $F_k$ as one ball soliton, whereas $B_{a_1,\cdots,a_r}$ is $r$ basket solitons.
For a composite basic soliton states, we consider the decomposition
in Proposition \ref{prop:composite} and sum over the disjoint union of pure solitons.  For a disjoint union of basic solitons, we sum over each basic soliton.
We also define amplitudes of the resulting pure solitons as follows;
amplitude of $F_k$ is $k$ and that of $B_a$ is $a$.

\begin{remark}
Alternatively, we may define the number of solitons of an arbitrary state
by the total number of pure solitons obtained by making enough many scatterings
of the state with ball solitons $F_k$ ($k$: large enough).
\end{remark}

\begin{example}
Scatterings of more than 60 $F_k$ ($k\geq 2$) with
$U_{10}B_7B_8U_{12}U_9FB_9F$ under $T_\ell$ ($\ell\geq 2$)
will give
$F\sqcup F\sqcup F\sqcup F\sqcup F\sqcup
B_{10,7,29,9}$.
Here the original composite state is cut into chunks as
$(U_{10})(B_7)(B_8U_{12}U_9F)(B_9)(F)$
and the above proposition applies to each chunk.
Thus the original state contains 5 ball solitons and 4 basket solitons.
\end{example}

By a $n$-scattering of solitons we mean the (large) time evolution of a disjoint union
$\cdots A_1 \cdots A_2 \cdots \cdots \cdots A_n \cdots$ of $n$ solitions $A_1, A_2, \ldots,A_n$
such that the solitons are arranged in decreasing speed from left to right.

\begin{theorem}\label{th:main}
The box-basket-ball system is solitonic in the following sense:
\begin{enumerate}
\item from any initial state the system eventually evolves into disjoint union of basic solitons,
\item the number and amplitudes of ball solitons and basket solitons contained in the initial state (after decomposition into pure solitons) are preserved under the time evolutions,
\item the scattering of $n$ solitons is factorized into two body scatterings.
\end{enumerate}
\end{theorem}
\begin{proof}
Consider the time evolution of the unbasketing of the system.
By Theorem \ref{thm:bb} at time $+ \infty$ the state consists of separate
solitons arranged in the order of non-decreasing speed.
The baskets cannot move with speed bigger than one.
Therefore after long enough time period all solitons in the unbasketing
 with speed greater than one shall come from solitons $F_k$ in the original box-basket-ball system.
The unbasketed solitons of speed one come from some combination of states $V$, $F$, $B_a$ and $U_a$ avoiding subsequences $FF$ and $FU$.
Thus they come from a union of basic speed one solitons in the box-basket-ball system.

For (2), let us consider a state $S$ which we assume is a disjoint union of basic solitons, arranged in decreasing order of speed.
In order to count the number of solitons, we consider a state
$S':=(F_k\otimes V^{\otimes k})^{\otimes M}\otimes V^{\otimes L}\otimes S$
where $M$ is an integer as large as the decomposition in Proposition \ref{prop:composite},
and $L$ is a much larger integer than $M$.
Then for a sufficiently large integer $N$, we can count the number of solitons from $T_\infty^N(S')$.  We may assume that $k$ is large enough, and that the spacing in $S$ is large enough that there is no scattering between the original basic solitons in the calculation of $T_\infty^N(S')$.

Now let us consider the state $T_\ell^{K}(S')$, where $K$ is large enough that the scattering between basic solitons in $S$ is completed, while $\ell$ is chosen small enough (relative to $k$ and $L$) that the $F_k$ have not yet scattered with $S$.
The number of solitons in the scattering of $S$ can then be counted by considering $T_\infty^NT_\ell^{K}(S') = T_\ell^K T_\infty^N(S')$ (see Theorem \ref{quantum_integrability}).  Finally, $T_\ell^K T_\infty^N(S')$ is obtained from $T_\infty^N(S')$ by the scattering of some pure solitons.  But by Proposition \ref{prop:FS} this preserves the number of ball solitons and basket solitons. The preservation of amplitudes also follows from this argument.

The proof of (3) is exactly the same as that of \cite[Theorem 4.6]{FOY},
and uses Theorem \ref{quantum_integrability} and Proposition \ref{prop:slow}:
we first evolve with $T_2$, which restricts to scatterings between speed $k >1$ and speed 1 solitons,
then we evolve with $T_3$, and so on.
After evolving with $T_2$, all the speed 1 solitons have been overtaken,
so that when we evolve with $T_3$, we will restrict ourselves to scattering of speed $k > 2$ and speed 2 solitons, and so on.
\end{proof}

This characterization of the solitonic property is defined and proved
for the box-ball system with capacity one boxes in \cite{TNS}.
Here we follow the treatment of \cite{FOY}.
For more general box-ball systems with boxes of arbitrary capacities,
it is proved in \cite{S1}.

\section{Proof of fast-slow scattering for $T_\infty$}
\label{sec:FS}
Let us break up time evolution into 
\begin{enumerate}
\item[(A)] Move baskets
\item[(B)] Move balls
\item[(C)] Reconfigure balls and baskets at each site separately
\end{enumerate}
Each time step will correspond to applying (A),(B),(C) in that order.
So that ``at time $t$'' we will have just completed step (C).
To emphasize this we sometimes say ``at integral time $t$''.

We shall adhere to the following rules when performing the moves (B) and (C).
In move (B), if a fast ball is placed inside a special basket,
it is designated a slow ball and the slow ball originally paired with
the special basket is now called a fast ball.
After the switch, we assume that the special basket and the new slow ball
(in the special basket before (C)) are paired.
In moves (B) or (C), we place balls in boxes first, then lower numbered baskets first,
and in move (C) if we have a choice, we shall place the slow ball in the box.

\subsection{In the middle of a scattering}

We shall first assume that no initial slow ball is present.

\begin{lemma}\label{lem:interval}
There exist integers $t_0, t_f$ and integers $i_t,j_t$ for $t \in [t_0,t_f]$ so that 
\begin{enumerate}
\item
$1 = i_{t_0} \leq j_{t_0} < i_{t_0+1} \leq j_{t_0+1} < \cdots < i_{t_f} \leq j_{t_f} =b$.
\item
At time $t \in [t_0,t_f]$ the fast balls are located inside either boxes or the non-special baskets with numbers inside $[i_t,j_t]$.
\item
At times $t < t_0$, no interaction between the solitons has occurred.
At times $t > t_f$, the scattering is complete and the output fast soliton is of the same length as before.
\end{enumerate}
\end{lemma}

\begin{lemma}\label{lem:phase} \
\begin{enumerate}
\item
Special baskets are always empty at integral time $t$.
Thus they do not experience phase shift.
\item
Non-special baskets are occupied at exactly one integral time.
Thus they are phase shifted by $-1$.
\item
Non-initial slow balls experience a phase shift of $-1$.
\end{enumerate}
\end{lemma}

\begin{proof}
We shall establish Lemma \ref{lem:interval} and \ref{lem:phase} by induction. 

We say that a special basket is activated when
a fast ball is placed inside it during move (B).
We first remark that assuming Lemma \ref{lem:phase}(1), we can show that
at integral times a special basket is always one site
behind its paired slow ball before activation,
and at the same site as its paired slow ball after activation.
Since when it is activated, the paired slow ball is the same
site with the paired special basket after (A).
Then the slow ball is interchanged with the fast ball in the spacial basket
when activation occurs.
Thus the slow ball enters the special basket and will not be moved by
the rest of the move (B) so that it stays at the same site with the paired special basket.

We show Lemma \ref{lem:interval}(1,2) for the first time $t = t_0$.
At time $t = t_0-1+(A)$, any special baskets have been moved
to the same site as the paired slow ball.
At $t = t_0-1 + (A,B)$, this special basket may or may not be occupied,
but after move (C), it will be unoccupied.
(The original paired slow ball will be moved,
and during move (C) the special basket will be emptied.)
Here the move (B) has the following property;
once the first fast ball is moved, then no slow balls
will be moved until all fast balls are moved.
Assuming there is no initial slow ball,
all slow balls are paired with special baskets.
Then if a fast ball is put into a special basket,
it is immediately interchanged by the paired slow ball and the new fast ball
will be moved by the rest of (B).

Now suppose the description of Lemma \ref{lem:interval}(1,2) holds at time $t$
and consider the next step.
In move (A), all except the non-special baskets in
$[i_{t},j_{t}]$ will be moved towards the right.
Then after the move (A), non-activated special baskets
are at the same site as the paired slow ball.
The already activated special baskets are one site ahead of the paired slow ball.
Note that we may assume by possibly changing the $j_{t}$
by 1 that any special baskets numbered in $[1,j_{t}]$ have been activated.

In move (B), the balls will move according to the following four steps:
\begin{enumerate}
\item
All slow balls to the left of any fast balls
will move to the same site as their activated special baskets.
\item
The first $k$ fast balls, where $k$ is the
number of activated special baskets numbered in the region $[i_t,j_t]$
will move into activated special baskets and immediately renamed as slow balls.
\item
The remainder of the fast balls, including any originally
slow balls (now renamed as fast) which are paired with the activated special
baskets in the second step or with newly activated special baskets
in the third step
will begin to fill up baskets labeled $j_{t}+1,j_{t}+2,\ldots$ and so on.
\item
The rest of the balls, all of which are slow
will move one step forward into a box.
\end{enumerate}
Just as in the time $t_0$ we can show that no slow balls
are touched during steps 2 and 3.
Thus the above four steps give the complete classification of
the move (B) except for the following cases.
In the final step it is possible
for a slow ball to end up in a basket rather than box;
and in the earlier steps it is possible for balls to end up in boxes
rather than baskets.
These exceptions do not affect the result of
Lemma \ref{lem:interval} and \ref{lem:phase}.


At time (C), every occupied special basket which has been activated will now be emptied
-- the original paired slow ball
(now designated a fast ball) has left the site,
and the slow ball inside the special basket can be placed into the box.
This shows that Lemma \ref{lem:interval}(1,2) holds at all times.
Lemma \ref{lem:interval}(3) also holds since the number of fast balls remain constant,
and they will come out all consecutive.

Lemma \ref{lem:phase}(1,2,3) all follow from the above discussion.
For (3), one just notes that the slow ball is always one ahead (at integral time)
of the paired special basket before activation, and always in the same site after activation.
\end{proof}

\subsection{Initial slow ball}
The analysis of the system containing an initial slow ball
can be reduced to that of no initial slow ball case.
The behaviour is just the same as the (known) scattering of $F_k$ and $F_1$.

\begin{lemma}
Any initial slow ball is phase shifted by $-2$, but does not affect the rest of the scattering.
\end{lemma}
\begin{proof}
When the initial slow ball is first overtaken by a fast ball, say at time $t + (A,B)$, we rename the left-most fast ball a slow ball (and now consider the initial slow ball a fast ball).  At time $t+2$, the new slow ball will be in the location of the original slow ball at time $t$.  This gives the calculation of the phase shift.  The rest of the scattering then proceeds in the original manner.  
\end{proof}

\subsection{Phase shift of fast soliton}
\begin{lemma}
Suppose there are $a$ baskets and $b$ balls in the slow soliton.
Then after scattering, the fast soliton is phase shifted by $2b-a$.
\end{lemma}
\begin{proof}
Every time a fast ball is placed into a non-special basket,
the fast soliton is slowed by $1$ unit.
When a special basket is first activated, the fast ball becomes a slow ball, and a new fast ball replaces the old one at the same site.  Since the new fast ball is occupying a box, this does not change the speed of the soliton, except that the new fast ball has yet to move.  And when the new fast ball moves, the fast soliton is sped up by $1$ unit.  The next time the activated special basket is occupied,
there is no effect on the speed of the fast balls
(the special basket traps a fast ball, but then a new fast ball will
replace it).
\end{proof}

\bigskip
\noindent
{\bf Acknowledgments.}
T.L. is supported by NSF grant DMS-0901111, and by a Sloan Fellowship.
P.P. was supported by NSF grant DMS-0757165.
R.S. is partially supported by Grants-in-Aid for Scientific Research
No. 21740114 from JSPS.

\end{document}